\documentclass[10pt]{fic-l}
\usepackage[]{amsmath}
\usepackage[]{amssymb}
\usepackage{amsfonts}

\newcommand{\abs}[1]{{\left|#1\right|}}

\newcommand{\RP}{\mathop\mathsf{RP}}
\newcommand{\PR}{\mathop\mathsf{PR}}
\newcommand{\RPPR}{\mathop\mathsf{RPPR}}
\newcommand{\ANY}{\mathop\mathsf{ANY}}
\newcommand{\ORD}{\mathop\mathsf{ORD}}
\newcommand{\ord}{{\mathop{\mathrm{ord}}\nolimits}}

\newcommand{\Nfp}{F}
\newcommand{\Mfp}{G}
\newcommand{\Ntc}{T}
\newcommand{\Nha}{C}
\newcommand{\sigzero}{d}
\newcommand{\sigone}{\sigma}

\renewcommand{\theenumi}{\alph{enumi}}
\renewcommand{\labelenumi}{(\theenumi)}
\renewcommand{\theenumii}{\roman{enumii}}
\renewcommand{\labelenumii}{\theenumii.}
\renewcommand{\theenumiii}{\Alph{enumiii}}
\renewcommand{\theenumiv}{\arabic{enumiv}}

\newtheorem{theorem}{Theorem}
\newtheorem{lemma}{Lemma}
\newtheorem{conjecture}{Conjecture}
\newtheorem{proposition}{Proposition}

{
\theoremstyle{remark}

}

\begin{document}

\title{New Conjectures and Results for Small Cycles of the Discrete Logarithm}
\author{Joshua Holden}
\thanks{The first author would like to thank the Rose-Hulman Institute of
Technology for the special stipend which supported this project.}
\address{Department of Mathematics,
Rose-Hulman Institute of Technology,
Terre Haute, IN, 47803-3999, USA}
\email{holden@rose-hulman.edu}

\author{Pieter Moree}
\address{Korteweg-de Vries Institute,
 Plantage Muidergracht 24, 1018 TV Amsterdam, The Netherlands}
\email{moree@science.uva.nl}
\date{\today}

\begin{abstract}
Brizolis asked the question: does every prime $p$ have a pair $(g,h)$ such
that $h$ is a fixed point for the discrete logarithm with base $g$?  The
first author previously extended this question to ask about not only fixed
points but also two-cycles, and gave heuristics (building on work of Zhang,
Cobeli, Zaharescu, Campbell, and Pomerance) for estimating the number of
such pairs given certain conditions on $g$ and $h$.  In this paper we give
a summary of conjectures and results which follow from these heuristics,
building again on the aforementioned work.  We also make some new
conjectures and prove some average versions of the results.
\end{abstract}

\maketitle

\section{Introduction and Statement of the Basic Equations}

Paragraph~F9  of~\cite{UPINT} includes the following problem,
attributed to Brizolis: given a prime $p>3$, is there always a
pair $(g,h)$ such that $g$ is a primitive root of $p$, $1 \leq h
\leq p-1$, and
\begin{equation} \label{fp}
g^{h} \equiv h  \mod{p} \enspace ?
\end{equation}
In other words, is there always a primitive root $g$ such that the
discrete logarithm $\log_{g}$ has a fixed point?  As we shall see,
Zhang~(\cite{Zhang}) not only answered the question for sufficiently
large $p$, but also estimated the number $N(p)$ of pairs $(g,h)$ which
satisfy the equation, have $g$ is primitive root, and also have $h$ a
primitive root which thus must be relatively prime to $p-1$.  This
result seems to have been discovered and proved by Zhang
in~\cite{Zhang} and later, independently, by Cobeli and Zaharescu
in~\cite{CZ}.  Campbell and Pomerance (\cite{CampbellTalk},
\cite{PomeranceTalk}) made the value of ``sufficiently large'' small
enough that they were able to use a direct search to affirmatively
answer Brizolis' original question.  As in~\cite{Holden02}, we will
also consider a number of variations involving side conditions on $g$
and $h$.

In~\cite{Holden02}, the first author also investigated the two-cycles of
$\log_{g}$, that is the pairs $(g,h)$ such that there is some $a$
between $1$ and $p-1$ such that
\begin{equation}\label{tc}
g^{h} \equiv a \mod{p} \quad \text{and} \quad g^{a} \equiv h \mod{p} .
\end{equation}
As we observed, attacking~(\ref{tc}) directly requires the
simultaneous solution of two modular equations, presenting both
computational and theoretical difficulties.  Whenever possible,
therefore, we instead work with the modular equation
\begin{equation} \label{ha}
    h^{h} \equiv a^{a} \mod{p} .
\end{equation}
Given $g$, $h$, and $a$ as in~(\ref{tc}), then (\ref{ha}) is clearly
satisfied and the common value is $g^{ah}$ modulo $p$.  Conditions on
$g$ and $h$ in~(\ref{tc}) can (sometimes) be translated into
conditions on $h$ and $a$ in~(\ref{ha}).  On the other hand, given a
pair $(h,a)$ which satisfies~(\ref{ha}), we can attempt to solve for
$g$ such that $(g,h)$ satisfies~(\ref{tc}) and translate conditions on
$(h,a)$ into conditions on $(g,h)$.  Again, we will investigate using
various side conditions.

Using the same notation as in~\cite{Holden02}, we will refer to an
integer which is a primitive root modulo $p$ as $\PR$ and an integer
which is relatively prime to $p-1$ as $\RP$.  An integer which is both
will be referred to as $\RPPR$ and one which has no restrictions will
be referred to as $\ANY$.  All integers will be taken to be between
$1$ and $p-1$, inclusive, unless stated otherwise.  If $N(p)$ is, as
above, the number of solutions to~(\ref{fp}) such that $g$ is a
primitive root and $h$ is a primitive root which is relatively prime
to $p-1$, then we will say $N(p)=\Nfp_{g \PR, h \RPPR}(p),$ and
similarly for other conditions.  Likewise the number of solutions
to~(\ref{tc}) will be denoted by $T$ and the number of solutions
to~(\ref{ha}) will be denoted by $C$.  If $\ord_{p}(g)=\ord_{p}(h)$,
we say that $g \ORD h$.

The idea of repeatedly applying the function $x \mapsto g^{x}
\bmod{p}$ is used in the famous cryptographically secure pseudorandom
bit generator of Blum and Micali.  (\cite{Blum-Micali}; see
also~\cite{Patel-Sundaram} and~\cite{Gennaro}, among others, for
further developments.)  If one could predict that a pseudorandom
generator was going to fall into a fixed point or cycle of small
length, this would obviously be detrimental to cryptographic security.
Our data suggests, however, that the chance that a pair $(g,h)$ is a
non-trivial two-cycle is $1/(p-1)$ for most of the conditions on
choosing $g$ and $h$ that we have investigated.  Likewise the chance
that a pair $(g,h)$ is a fixed point is generally $1/(p-1)$. This
might perhaps be taken as an indication that the seed of one of these
pseudorandom generators should be chosen to avoid redundant conditions 
which would increase the chances of a small cycle.

This paper is meant to serve as a summary of the authors' recent 
work.  For detailed proofs and explanation we refer the reader to 
our forthcoming paper~(\cite{HM}), in preparation.  Numerical examples are 
provided here to illustrate the conjectures and results.

\section{Conjectures and Theorems for Fixed Points}

\label{fpconjsec}

A list of conjectures and theorems on fixed points appeared
in~\cite{Holden02} and was corrected in the unpublished
notes~\cite{Holden02a}.  These conjectures and theorems are summarized
in Table~\ref{fptalktable}, which appeared in~\cite{Holden02a}.  The
table also contains new data collected since~\cite{Holden02}.

\begin{table}[!ht]
    \caption{Solutions to~(\ref{fp})}
    \label{fptalktable}
    $$\begin{array}{|l|l|l|l|l|}
\multicolumn{5}{l}{\text{(a) Predicted formulas for $\Nfp(p)$}}\\
    \hline
        g \setminus h & \ANY & \PR & \RP & \RPPR  \\
        \hline
    \ANY &\approx \scriptstyle(p-1) &
    \approx\frac{\phi(p-1)^{2}}{(p-1)} &
        = \scriptstyle \phi(p-1) &
    \approx\frac{\phi(p-1)^{2}}{(p-1)} \\

        \hline
    \PR & \approx \scriptstyle \phi(p-1) &
    \approx\frac{\phi(p-1)^{2}}{(p-1)} &
    \approx\frac{\phi(p-1)^{2}}{(p-1)} &
    \approx\frac{\phi(p-1)^{2}}{(p-1)} \\

    \hline
    \RP & \approx \scriptstyle\phi(p-1) &
    \approx\frac{\phi(p-1)^{3}}{(p-1)^2} &
    \approx\frac{\phi(p-1)^{2}}{(p-1)} &
    \approx\frac{\phi(p-1)^{3}}{(p-1)^2} \\

    \hline
    \RPPR & \approx\frac{\phi(p-1)^{2}}{(p-1)} &
    \approx\frac{\phi(p-1)^{3}}{(p-1)^2} &
    \approx\frac{\phi(p-1)^{3}}{(p-1)^2} &
    \approx\frac{\phi(p-1)^{3}}{(p-1)^2} \\
        \hline
\multicolumn{5}{l}{}\\
    \multicolumn{5}{l}{\text{(b) Predicted values for
    $\Nfp(100057)$}}\\
        \hline
        g \setminus h & \ANY & \PR & \RP & \RPPR  \\
        \hline
\ANY & 100056 & 9139.46 & 30240 & 9139.46 \\
    \hline
    \PR  &30240  &9139.46  &9139.46 & 9139.46 \\
    \hline
    \RP  &30240  &2762.23 & 9139.46 & 2762.23 \\
    \hline
    \RPPR  &9139.46&  2762.23 & 2762.23 & 2762.23  \\
    \hline
\multicolumn{5}{l}{}\\
    \multicolumn{5}{l}{\text{(c) Observed values for
    $\Nfp(100057)$}}\\

        \hline
        g \setminus h & \ANY & \PR & \RP & \RPPR  \\
        \hline
   \ANY&  98506& 9192 & 30240& 9192\\
   \hline
\PR&   29630& 9192&  9192&  9192\\
\hline
\RP&   29774& 2784&  9037&  2784\\
\hline
\RPPR& 9085&  2784&  2784&  2784\\
\hline
\end{array}$$
\end{table}

\label{fpthmsec}

The first rigorous result on this subject was for $\Nfp_{g
\PR, h \RPPR}(p)$.  Both~\cite{Zhang} and~\cite{CZ} provided bounds on
the error involved; we will use notation closer to~\cite{CZ}.

\begin{theorem}[Theorem~1 of~\cite{CZ}] \label{cz1thm}
\[
\abs{\Nfp_{g \PR, h \RPPR}(p)  - \frac{\phi(p-1)^{2}}{p-1}} \leq
\sigzero(p-1)^{2}\sqrt{p}(1+\ln p).
\]
\end{theorem}

We next turn our attention to $\Nfp_{g \ANY, h \ANY}(p)$, for 
which we can prove the following result:

\begin{theorem} \label{thm3}
    \[
\abs{\Nfp_{g \ANY, h \ANY}(p)  - (p-1)} \leq
 \sigzero(p-1) \sigone(p-1) \sqrt{p}(1+\ln p).
\]
\end{theorem}

Unfortunately, $\sigone(n) = O(n \ln \ln n)$ in the worst case and in
any case $\sigone(p-1) \geq p-1 + (p-1)/2 + 2 + 1 > 3p/2$.  Thus the
error term overwhelms the main term.  The problem occurs because we
use the fact that~(\ref{fp}) can be solved exactly when $\gcd(h,
p-1)=e$ and $h$ is a $e$-th power modulo $p$, and in fact there are
exactly $e$ such solutions.  When $h$ is $\RPPR$ then $e$ is always
$1$ so counting the number of $h$ is sufficient.  When $h$ is $\ANY$,
however, we need to count the number of $h$ such that $\gcd(h, p-1)=e$
and $h$ is a $e$-th power modulo $p$ and then multiply by $e$, and do
this for each divisor $e$ of $p-1$.  Thus an error of even $1$ in
calculating the number of $h$ above for a large value of $e$ will
result in an error of $O(p-1)$.  (We can improve the situation
somewhat by separating out the elements where $e$ is $p-1$ or
$(p-1)/2$, but the results are still not what one would wish for.  
More details will appear in~\cite{HM}.)

%

%
%


The case where $g$ is $\PR$ and $h$ is $\ANY$  is very
similar to the previous case, and unfortunately has the same problem:

\begin{theorem} \label{thm4}
       \[
\abs{\Nfp_{g \PR, h \ANY}(p)  - \phi(p-1)} \leq
 \sigzero(p-1)^{2} \sigone(p-1) \sqrt{p}(1+\ln p).
\]
\end{theorem}

Finally, we should mention that the second author (in~\cite{MoreeMR})
pointed out that we could also estimate the number $\Mfp_{g
\PR, h \ANY}(p)$ of values $h$ such that there exists \emph{some} $g$
satisfying~(\ref{fp}), with $g \PR$ and $h \ANY$:

\begin{theorem} \label{thm5}
    \[
\abs{\Mfp_{g \PR, h \ANY}(p)  - \frac{1}{p-1} \sum_{e \mid
p-1} \phi \left( \frac{p-1}{e} \right)^{2}} \leq
\sigzero(p-1)^{3}\sqrt{p}(1+\ln p).
\]
\end{theorem}

Similarly, we have:

\begin{theorem} \label{thm6}
    \[
\abs{\Mfp_{g \ANY, h \ANY}(p)  - \sum_{e \mid
p-1}\frac{1}{e}  \phi \left( \frac{p-1}{e} \right)} \leq
\sigzero(p-1)^{2}\sqrt{p}(1+\ln p).
\]
\end{theorem}

Since we are no longer counting multiple solutions for each value 
of $h$ the problem mentioned above disappears; the error terms are 
$O(p^{1/2+\epsilon})$ while the main terms look on average like a 
constant times $p$.


\section{Conjectures for Two-Cycles}

Conjectures relating to equations~(\ref{ha}) and~(\ref{tc}) also
appeared in~\cite{Holden02} and were corrected in the unpublished
notes~\cite{Holden02a}.  These are summarized in
Tables~\ref{hatalktable} and~\ref{tctalktable}, which appeared
in~\cite{Holden02a}.  The table also contains new data collected
since~\cite{Holden02}.  As in~\cite{Holden02}, we distinguish between
the ``trivial'' solutions to~(\ref{ha}), where $h=a$, and the
``nontrivial'' solutions.

It was observed in~\cite{Holden02} that when neither $h$ nor $a$ is $\RP$
the relationship between~(\ref{tc}) and~(\ref{ha}) is more complicated than
in the other cases.  (Summaries of the conjectures in these cases are given
in Tables~\ref{hatalktable} and~\ref{tctalktable}.)  We were able, however,
to make the following conjectures about solutions to~(\ref{ha}).

\begin{conjecture} \label{conj7} \mbox{}
    \begin{enumerate}
\item \label{conj7a} $\Nha_{h \ANY, a \ANY}(p) \approx (p-1) + \sum_{m
\mid p-1} {\phi(m)} \left( \sum_{d \mid (p-1)/m}
\frac{\phi(dm)}{dm}\right)^{2}$.

\item \label{conj7b} If $p-1$ is squarefree then $\Nha_{h \ANY, a
\ANY}(p) \approx (p-1) + \prod_{q \mid p-1} \left( q + 1 - \frac{1}{q}
\right)$, where the product is taken over primes $q$ dividing $p-1$. 

\item \label{conj7c} In general,
    \begin{multline*} \label{geneq}
    \Nha_{h \ANY, a \ANY}(p)  \\
\shoveleft{\approx (p-1) +  \prod_{q^\alpha \Vert p-1} \left(
\left[\left(1-\frac{1}{q}\right)\alpha + 1\right]^2 \right.} \\
+ \left(1-\frac{1}{q}\right)^3
\left[ (\alpha+1)^2 \frac{q^{\alpha+1}-q}{q-1}
- 2 (\alpha+1) \frac{\alpha q^{\alpha+2} -
(\alpha+1)q^{\alpha+1}+q}{(q-1)^2} \right.\\
\left.\left.
+ \frac{\alpha^2 q^{\alpha+3} - (2\alpha^2 + 2\alpha-1)q^{\alpha+2}
        + (\alpha^2+2\alpha+1)q^{\alpha+1}-q^2-q}{(q-1)^3}
\right]
\right),
\end{multline*}
where the product is taken over primes $q$ dividing $p-1$ and $\alpha$ is
the exact power of $q$ dividing $p-1$.

\item \label{conj7d} $\Nha_{h \PR, a \ANY}(p) \approx 2 \phi(p-1)$.
\item \label{conj7e}
$\Nha_{h \ANY, a \PR}(p)\approx 2 \phi(p-1)$.
\item \label{conj7f}
 $\Nha_{h \PR, a \PR}(p) \approx \phi(p-1) +
\phi(p-1)^{2}/(p-1)$.

\end{enumerate}
\end{conjecture}

(The formulas in Conjecture~\ref{conj7}(\ref{conj7a}) and
Conjecture~\ref{conj7}(\ref{conj7c})
appear in~\cite{Holden02} with typos.  They appear correctly here and
in~\cite{Holden02a}.)

\begin{table}[!ht]
    \caption{Solutions to~(\ref{ha})}
    \label{hatalktable}

    $$\begin{array}{|l|l|l|l|l|}

\multicolumn{5}{l}{\text{(a) Predicted formulas for the nontrivial
part of $\Nha(p)$}}\\
    \hline
        a \setminus h & \ANY & \PR & \RP & \RPPR  \\
        \hline
    \ANY & \approx 
\sum {\scriptstyle\phi(m)} \left( \sum
\frac{\phi(dm)}{dm}\right)^{2} &
\approx \scriptstyle \phi(p-1) &
    \approx \scriptstyle \phi(p-1) &
    \approx \frac{\phi(p-1)^{3}}{(p-1)^{2}} \\

        \hline
    \PR & \approx \scriptstyle \phi(p-1) &
    \approx \frac{\phi(p-1)^{2}}{(p-1)} &
    \approx \frac{\phi(p-1)^{2}}{(p-1)} &
    \approx \frac{\phi(p-1)^{3}}{(p-1)^{2}} \\

    \hline
    \RP &  \approx \scriptstyle \phi(p-1) &
    \approx \frac{\phi(p-1)^{2}}{(p-1)} &
    \approx \frac{\phi(p-1)^{2}}{(p-1)} &
    \approx \frac{\phi(p-1)^{3}}{(p-1)^{2}} \\

    \hline
    \RPPR & \approx \frac{\phi(p-1)^{3}}{(p-1)^{2}} &
    \approx \frac{\phi(p-1)^{3}}{(p-1)^{2}} &
    \approx \frac{\phi(p-1)^{3}}{(p-1)^{2}} &
    \approx \frac{\phi(p-1)^{3}}{(p-1)^{2}} \\

        \hline
\multicolumn{5}{l}{}\\
\multicolumn{5}{l}{\text{(b) Predicted values for the nontrivial
part of $\Nha(100057)$}}\\
        \hline
        a \setminus h & \ANY & \PR & \RP & \RPPR  \\
        \hline
\ANY & 190822.0 & 30240 & 30240 & 2762.225  \\
\hline
\PR & 30240 & 9139.458 & 9139.458 & 2762.225  \\
\hline
\RP & 30240 & 9139.458 & 9139.458 & 2762.225  \\
\hline
\RPPR & 2762.225 & 2762.225 & 2762.225 & 2762.225 \\
    \hline
\multicolumn{5}{l}{}\\
\multicolumn{5}{l}{\text{(c) Observed values for the nontrivial
part of $\Nha(100057)$}}\\
        \hline
        a \setminus h & \ANY & \PR & \RP & \RPPR  \\
        \hline
\ANY &  190526 & 30226 & 30291 & 2820 \\
\hline
\PR &   30226 & 9250 & 9231 & 2820 \\
\hline
\RP &   30291 & 9231 & 9086 & 2820\\
\hline
\RPPR & 2820 & 2820 & 2820 & 2820 \\
\hline
\end{array}$$
\end{table}

As observed in~\cite{Holden02}, conditions on~(\ref{tc}) can sometimes
be translated into conditions on~(\ref{ha}) in a relatively
straightforward manner.  In other cases, however, things are more
complicated. 
Let $d=\gcd(h,a,p-1)$, and let $u_{0}$ and $v_{0}$ be such that
\[
u_{0} h + v_{0} a \equiv d \mod{p-1}.
\]
Taking the logarithm of the two equations of~(\ref{tc}) with respect
to the same primitive root $b$ and using Smith Normal Form, we can
show that~(\ref{tc}) is equivalent to the equations:

\begin{equation} \label{tcsmith}
    h^{h/d} \equiv a^{a/d} \mod{p} 
    \quad \text{and} \quad
    g^{d}\equiv h^{v_{0}} a^{u_{0}} \mod{p}.
\end{equation}
In the case where $d=\gcd(h,a,p-1)=1$ then this becomes just

\begin{equation}
    h^{h} \equiv a^{a} \mod{p} 
    \quad \text{and} \quad
    g  \equiv h^{v_{0}} a^{u_{0}} \mod{p}.
\end{equation}

Thus:

\begin{proposition} \label{equivprop}
    If $\gcd(h,a,p-1)=1$, then there is a one-to-one correspondence
    between triples $(g,h,a)$ which satisfy~{\rm{(\ref{tc})}} and pairs
    $(h,a)$ which satisfy~{\rm{(\ref{ha})}}, and the value of $g$ is unique
    given $h$ and $a$.  In particular, this is true if $h$ is $\RP$ or
    $a$ is $\RP$.
\end{proposition}

In~\cite{Holden02} it was claimed that given a pair $(h,a)$ which is a
solution to~(\ref{ha}) we expect on the average
$\gcd(a,p-1)\gcd(h,p-1)/\gcd(ha, p-1)^{2}$ pairs $(g,h)$ which are
solutions to~(\ref{tc}).
It is clear from~(\ref{tcsmith}), however, that the proper equation to
look at in this case is not~(\ref{ha}), but
\begin{equation} \label{hasmith}
    h^{h/d} \equiv a^{a/d} \mod{p}.
\end{equation}

Now we can approximate the number of nontrivial solutions
of~(\ref{hasmith}) using a similar birthday paradox argument to that
used in~\cite{Holden02} for Conjecture~\ref{conj7}.  The end result 
(see our forthcoming paper for details) is
the following conjectures:
\begin{conjecture} \label{conj8} \mbox{}
    \begin{enumerate}
    \item
$\Ntc_{g \PR, h \ANY}(p) \approx 2\phi(p-1)$.
\item
    $\Ntc_{g \ANY, h \ANY}(p) \approx
2(p-1)$.
\end{enumerate}
\end{conjecture}

and:

\begin{conjecture} \label{conj9} \mbox{} 
\begin{enumerate}
    \item $\Ntc_{g \RP, h \bullet}(p) \approx \left[\phi(p-1)/(p-1)\right]
    \Ntc_{g \ANY, h \bullet}(p)$.  
    
    \item $\Ntc_{g \RPPR, h \bullet}(p) \approx
    \left[\phi(p-1)/(p-1) \right] \Ntc_{g \PR, h \bullet}(p)$.
\end{enumerate}
(where $\bullet$ stands for any one of the four conditions which we have 
used on $h$)
\end{conjecture}

The data from Tables~\ref{fptalktable}, \ref{hatalktable},
and~\ref{tctalktable} was collected on a Beowulf cluster\footnote{A
type of high-speed parallel computing system built out of standard PC
parts.}, with 19 nodes, each consisting of 2 Pentium~III processors
running at 1~Ghz.  The programming was done in C, using MPI, OpenMP,
and OpenSSL libraries.  The collection took 68 hours for all values of
$\Nfp(p)$, $\Ntc(p)$, and $\Nha(p)$, for five primes $p$ starting at
100000.

\begin{table}[!ht]
    \caption{Solutions to~(\ref{tc})}
    \label{tctalktable}

   $$\begin{array}{|l|l|l|l|l|}
\multicolumn{5}{l}{\text{(a) Predicted formulas for the nontrivial
part of $\Ntc(p)$}}\\
    \hline
        g \setminus h & \ANY & \PR & \RP & \RPPR  \\
        \hline
    \ANY & \approx \scriptstyle (p-1)    &
    \approx \frac{\phi(p-1)^{2}}{(p-1)} &
    \approx \scriptstyle \phi(p-1)  &
    \approx \frac{\phi(p-1)^{3}}{(p-1)^{2}} \\

        \hline
    \PR & \approx \scriptstyle \phi(p-1)   &
     \approx \frac{\phi(p-1)^{2}}{(p-1)} &
    \approx \frac{\phi(p-1)^{2}}{(p-1)} &
   \approx \frac{\phi(p-1)^{3}}{(p-1)^{2}} \\

    \hline
    \RP &  \approx \scriptstyle \phi(p-1) &
    \approx \frac{\phi(p-1)^{3}}{(p-1)^{2}} &
     \approx \frac{\phi(p-1)^{2}}{(p-1)} &
   \approx \frac{\phi(p-1)^{4}}{(p-1)^{3}} \\

    \hline
    \RPPR & \approx \frac{\phi(p-1)^{2}}{(p-1)}     &
    \approx \frac{\phi(p-1)^{3}}{(p-1)^{2}} &
    \approx \frac{\phi(p-1)^{3}}{(p-1)^{2}} &
    \approx \frac{\phi(p-1)^{4}}{(p-1)^{3}} \\

        \hline
\multicolumn{5}{l}{}\\
\multicolumn{5}{l}{\text{(b) Predicted values for the nontrivial
part of $\Ntc(100057)$}}\\
        \hline
        g \setminus h & \ANY & \PR & \RP & \RPPR  \\
        \hline
\ANY & 100056 & 9139.5 & 30240 & 2762.2  \\
\hline
\PR &30240 & 9139.5 & 9139.5 & 2762.2 \\
\hline
\RP & 30240 & 2762.2 & 9139.5 & 834.8   \\
\hline
\RPPR & 9139.5 & 2762.2 & 2762.2 & 834.8  \\
    \hline
\multicolumn{5}{l}{}\\
\multicolumn{5}{l}{\text{(c) Observed values for the nontrivial
part of $\Ntc(100057)$}}\\
        \hline
        g \setminus h & \ANY & \PR & \RP & \RPPR  \\
        \hline
\ANY & 100860 & 9231 & 30291 & 2820\\
\hline
\PR & 30850 & 9231 & 9231 & 2820\\
\hline
\RP &    30368 & 2882 & 9240 & 916\\
\hline
\RPPR & 9376 & 2882 & 2882 & 916\\

    \hline
\end{array}$$
\end{table}

\section{Averages of the Results and Conjectures}
 
Thus far we have considered variants of Brizolis conjecture for a
fixed finite field with $p$ elements.  In this section we consider
average versions of these results and conjectures.  The conjectures
predict a main term; the results give a main term and an error term. 
The following sequence of lemmas gives the behavior of the main terms,
on average.
 
The following result for $k=1$ is well-known, see e.g. \cite{MoreeJNT,
Stephens}.  For arbitrary $k$ it was claimed by Esseen \cite{Esseen}
(but only proved for $k=3$).  A proof can be given based on an idea of
Carl Pomerance \cite{PomerancePC}.  (Proofs of all of the results in
this section will appear in a forthcoming paper.)

\begin{lemma}
    \label{kinderlijk1}
    Let $k$ and $C$ be arbitrary real numbers with $C>0$. Then
    $$ 
    \sum_{p\le x}\left(\phi(p-1)\over p-1\right)^k=A_k~{\rm Li}(x)
    +O_{C,k}\left({x\over \log^C x}\right),
    $$ 
    where $$A_k=\prod_p\left(1+{(1-1/p)^k-1\over p-1}\right).$$
\end{lemma}

Given this lemma it is trivial to establish:

\begin{theorem} Let $C>0$ be arbitrary.  We have $$\sum_{p\le
x}{\frac{\Nfp_{g \PR,h \RPPR}(p)}{p-1}}=A_{2}{\rm
Li}(x)+O_C\left({x\over \log^Cx}\right).$$
\end{theorem}

Using similar lemmas, one can prove:

\begin{theorem} Let $C>0$ be arbitrary.  We have
$$\sum_{p\le
x}{\Mfp_{g{\rm PR},h{\rm ANY}}(p)\over p-1}=A_1{\zeta(3)\over
\zeta(2)}{\rm
Li}(x)+O_C\left({x\over \log^Cx}\right),$$
where $$A_1{\zeta(3)\over \zeta(2)}=\prod_p\left(1-{2p\over
p^3-1}\right)\approx 0.27327~30607~85299~15983\cdots$$
and $$\sum_{p\le x}{\Mfp_{g{\rm ANY},h{\rm ANY}}(p)\over p-1}=S{\rm
Li}(x)+O_C\left({x\over \log^Cx}\right),$$
where $$S=\prod_p\left(1-{p\over p^3-1}\right)\approx
0.57595~99688~92945~43964\cdots $$
is the Stephens constant (see~\cite{Stephens76}).

\end{theorem}

Theorems~\ref{thm3} and~\ref{thm4} are 
unfortunately more problematic, due to the presence of the 
exceptionally large error term.  The error term can 
probably be reduced to no larger order than the main term by 
separating out the most problematic cases and considering the 
sort of averaging we are doing in this section but the results are 
still conjectural at present, and the error term is still not 
satisfactory in any case.

On the other hand, almost all of the conjectures on~(\ref{fp}), (\ref{ha}),
and~(\ref{tc}) lend themselves easily to average versions of the sort
treated above.  These average versions are summarized in
Tables~\ref{fpavetable}, \ref{haavetable}, and~\ref{tcavetable}.  The data
in these tables was collected on the same Beowulf cluster mentioned above,
with similar software.  The collection took 17 hours for all values of
$\sum_{p\leq x}\frac{\Nfp(p)}{p-1}$, $\sum_{p\leq x}\frac{\Ntc(p)}{p-1}$,
and $\sum_{p\leq x}\frac{\Nha(p)}{p-1}$, for $x=6143$.

The results of the preceding section unfortunately do not allow
us to evaluate the average value of the right hand side of 
Conjecture~\ref{conj7}(\ref{conj7a}).
Let us put
$$w(p)=\sum_{m|p-1}\phi(m)\left(\sum_{d|m}{\phi(dm)\over
dm}\right)^2.$$
Numerically it seems that
$$\lim_{x\rightarrow \infty}{1\over \pi(x)}\sum_{p \leq x}{w(p)\over p-1}
=1.644\cdots,$$
with rather fast convergence. We are thus tempted to propose the
following conjecture.

\begin{conjecture}
 Let $C>0$ be arbitrary.  We have 
 
 $$\sum_{p\le x}{\Nha_{h \ANY, a \ANY}(p)\over p-1}=2.644\cdots{\rm
 Li}(x)+O_C\left({x\over \log^Cx}\right).$$
\end{conjecture}

Although we cannot prove this at present, we can establish the following result.
\begin{lemma}
For every $x$ sufficiently large we have
$$1.444\le {1\over 
\pi(x)}\sum_{p \leq x}{w(p)\over p-1}\le 3.422$$
\end{lemma}

\begin{table}[!ht]
    \caption{Average Solutions to~(\ref{fp})}
    \label{fpavetable}
    $$\begin{array}{|l|l|l|l|l|}
\multicolumn{5}{l}{\text{(a) Predicted approximate values for $\frac{1}{\pi(x)}
\sum_{p \leq x} \Nfp(p)$}}\\
    \hline
        g \setminus h & \ANY & \PR & \RP & \RPPR  \\
        \hline
	\ANY & 1 & 0.1473494000 & 0.3739558136 & 0.1473494000 \\

        \hline
    \PR &0.3739558136 & 0.1473494000 &0.1473494000 & 0.1473494000 \\

    \hline
    \RP & 0.3739558136 & 0.0608216551 & 0.1473494000 & 0.0608216551 \\

    \hline
    \RPPR &0.1473494000 & 0.0608216551 & 0.0608216551 & 0.0608216551\\
        \hline
\multicolumn{5}{l}{}\\
    \multicolumn{5}{l}{\text{(b) Observed values for
    $x=6143$}}\\

        \hline
        g \setminus h & \ANY & \PR & \RP & \RPPR  \\
        \hline
   \ANY&  0.9904034375&0.14851987375&0.37592474125&0.14851987375\\

   \hline
\PR&   0.3749536975&0.14851987375&0.14851987375&0.14851987375\\
\hline
\RP&   0.3739629175&0.0612404775&0.15122619375&0.0612404775\\
\hline
\RPPR& 0.14792889125&0.0612404775&0.0612404775&0.0612404775\\

\hline
\end{array}$$
\end{table}

\begin{table}[!ht]
    \caption{Average Solutions to~(\ref{ha})}
    \label{haavetable}

    $$\begin{array}{|l|l|l|l|l|}

\multicolumn{5}{l}{\text{(a) Predicted approximate values for the nontrivial
part of $\frac{1}{\pi(x)}
\sum_{p \leq x} \Nha(p)$}}\\  
    \hline
        a \setminus h & \ANY & \PR & \RP & \RPPR  \\
        \hline
    \ANY & 1.644\cdots &0.3739558136 & 0.3739558136 &
    0.0608216551 \\

        \hline
    \PR &  0.3739558136 & 0.1473494000& 0.1473494000 &  0.0608216551 \\

    \hline
    \RP & 0.3739558136 & 0.1473494000 &  0.1473494000 & 0.0608216551 \\

    \hline
      \RPPR & 0.0608216551 & 0.0608216551 & 0.0608216551& 0.0608216551\\

    \hline
\multicolumn{5}{l}{}\\
\multicolumn{5}{l}{\text{(b) Observed values for the nontrivial part 
for
    $x=6143$}}\\
        \hline
        a \setminus h & \ANY & \PR & \RP & \RPPR  \\
        \hline
\ANY & 1.6113896337 & 0.3655877485 & 0.3765792535 & 0.060552674 \\
\hline
\PR & 0.3655877485 & 0.14608992975 & 0.1478925015 & 0.060552674 \\
\hline
\RP & 0.3765792535 & 0.1478925015 & 0.146740421 & 0.060552674 \\
\hline
\RPPR &0.060552674 & 0.060552674 & 0.060552674 & 0.060552674\\
\hline
\end{array}$$
\end{table}

\begin{table}[!ht]
    \caption{Average Solutions to~(\ref{tc})}
    \label{tcavetable}

   $$\begin{array}{|l|l|l|l|l|}
\multicolumn{5}{l}{\text{(a) Predicted approximate values for the nontrivial
part of $\frac{1}{\pi(x)}
\sum_{p \leq x} \Ntc(p)$}}\\  
    \hline
        g \setminus h & \ANY & \PR & \RP & \RPPR  \\
        \hline
    \ANY & 1 & 0.1473494000 & 0.3739558136 & 0.0608216551\\

        \hline
    \PR &0.3739558136 & 0.1473494000 &0.1473494000 & 0.0608216551 \\

    \hline
    \RP & 0.3739558136 & 0.0608216551 &0.1473494000 & 0.0261074463\\

    \hline
    \RPPR & 0.1473494000 & 0.0608216551 & 0.0608216551 & 0.0261074463 \\

        \hline
\multicolumn{5}{l}{}\\
\multicolumn{5}{l}{\text{(b) Observed values for the nontrivial part 
for
    $x=6143$}}\\
        \hline
        g \setminus h & \ANY & \PR & \RP & \RPPR  \\
        \hline
\ANY & 0.9933146575 & 0.14884923375 & 0.3772284725 & 0.06150940625\\
\hline
\PR & 0.37381320625 & 0.14884923375 & 0.14884923375 & 0.06150940625\\
\hline
\RP & 0.36701980375 & 0.06089004625 & 0.146029115 & 0.02640389625\\
\hline
\RPPR & 0.14697618875 & 0.06089004625 & 0.06089004625 & 0.02640389625\\

    \hline
\end{array}$$
\end{table}

\section*{Acknowledgments}

Once again, the first author would like to thank the people mentioned
in~\cite{Holden02}: John Rickert, Igor Shparlinski, Mariana Campbell,
and Carl Pomerance.  He would also like to thank Victor Miller
for the suggestion to use the Smith Normal Form.

The authors would like to thank the anonymous referees for several 
helpful suggestions.

\end{document}